\documentclass[aoas,MSNbibl,nameyear,rotating,dvips]{arximspdf}
\usepackage{multirow,mathbh,dcolumn}
\usepackage{graphicx}

\doi{10.1214/13-AOAS638} 
\volume{7}
\issue{3}
\pubyear{2013}
\firstpage{1733}
\lastpage{1762}

\makeatletter
 \newcolumntype{d}[1]{D{.}{.}{#1}}
\newcommand{\eqref}[1]{(\ref{#1})}

\newtheorem{Algorithm}{Algorithm}

\renewcommand{\epsilon}{\varepsilon}
\newcommand{\R}{\mathbb{R}}
\renewcommand{\vec}[1]{\mathbf{#1}}
\makeatother

\begin{document}
\begin{frontmatter}

\title{A method for generating realistic correlation matrices}
\runtitle{Generating correlation matrices}

\begin{aug}
\author[a]{\fnms{Johanna} \snm{Hardin}\corref{}\thanksref{t1}\ead[label=e1]{jo.hardin@pomona.edu}},
\author[a]{\fnms{Stephan Ramon} \snm{Garcia}\thanksref{t2}\ead[label=e2]{Stephan.Garcia@pomona.edu}}
\and\break
\author[b]{\fnms{David} \snm{Golan}\thanksref{t1,t3}\ead[label=e3]{davidgo5@post.tau.ac.il}}
\thankstext{t1}{Supported in part by the Institute for Pure and Applied
Mathematics, NSF Grant DMS-09-31852.}
\thankstext{t2}{Supported in part by NSF Grant DMS-10-01614.}
\thankstext{t3}{Supported in part by a fellowship from the Edmond J.
Safra center for Bioinformatics at Tel Aviv University.}

\runauthor{J. Hardin, S.~R. Garcia and D. Golan}
\affiliation{Pomona College, Pomona College and Tel Aviv University}

\address[a]{J. Hardin\\
S.~R. Garcia\\
Department of Mathematics\\
Pomona College\\
610 N. College Ave.\\
Claremont, California 91711\\
USA\\
\printead{e1}\\
\phantom{E-mail:\ }\printead*{e2}}

\address[b]{D. Golan\\
Department of Statistics and OR\\
School of Mathematical Sciences\\
Tel Aviv University\\
69975 Tel Aviv\\
Israel\\
\printead{e3}}
\end{aug}

\received{\smonth{8} \syear{2012}}
\revised{\smonth{1} \syear{2013}}

%
\begin{abstract}
Simulating sample correlation matrices is important in many areas of
statistics. Approaches such as generating Gaussian data and finding
their sample correlation matrix or generating random uniform $[-1,1]$
deviates as pairwise correlations both have drawbacks. We develop an
algorithm for adding noise, in a highly controlled manner, to general
correlation matrices. In many instances, our method yields results
which are superior to those obtained by simply simulating Gaussian
data. Moreover, we demonstrate how our general algorithm can be
tailored to a number of different correlation models. Using our results
with a few different applications, we show that simulating correlation
matrices can help assess statistical methodology.
\end{abstract}

%

\end{frontmatter}

\section{Introduction}\label{sec1}

As computational resources continue to improve, researchers can take
advantage of simulation studies to investigate properties and results
associated with novel statistical methodology. In particular,
simulating correlation matrices with or without a given structure can
provide insight into the sensitivity of a model. There has been
extensive work on simulating correlation matrices with random entries,
that is, generating positive-semidefinite matrices with all entries
bounded by $[-1,1]$ and having ones along the diagonal. Seminal work by
\citet{Marsaglia84} discusses distributional characteristics and
eigenvalues of simulated random correlation matrices. Although there
has been additional work expanding the ideas associated with generating
random correlation matrices [\citet{Joe06},
\citet{Lewandowski09},
\citet{Holmes91},
\citet{Davies00},
\citet{Rae97}] and even randomly
generating correlation matrices within particular settings [\citet
{Ng10},
\citet{Holmes89}], to our knowledge there is no literature devoted to
the problem of adding noise to given template correlation structures.

We discuss the need to simulate realistic correlation matrices in a
specific context. By \emph{realistic} we mean not only that the
correlation matrix has some prescribed structure (dependent upon the
requirements of the particular application), but also that it is \emph
{noisy}. Below, we discuss the importance of simulating correlation
matrices in probit analysis, Bayesian shrinkage estimation,
meta-analysis, multiple comparisons, management science, factor
analysis, heritability estimation, network analysis and classification.

In order to ensure identifiability of model parameters, covariance
matrices in probit analysis on longitudinal data are often constrained
to be correlation matrices. \citet{Liu06} and \citet{Zhang06} discuss
advantages and disadvantages of different prior distributions used with
a Metropolis Hastings algorithm to sample correlation matrices from a
posterior distribution.

\citet{Barnard00} use hierarchical Bayesian models to describe the
covariance between parameters in a regression model with shrinkage. In
practice, joint uniform and inverse-Wishart priors are used to simulate
correlation matrices.

One important task in meta-analysis is to combine correlation matrices
from different studies. Different methods for combining such matrices
are used to make inferences about marginal correlations. \citet
{Hafdahl07} ran a Monte Carlo study generating sample correlation
matrices using Gaussian deviates from a given population correlation matrix.

One of the big challenges in developing accurate multiple comparisons
techniques is knowing the underlying correlation structure of the many
items being compared. Simply knowing the rate of null versus
alternative tests in a given simulation does not provide enough
information for realistic application to studies with possibly strongly
correlated data and hypotheses. In order to model false discovery rates
(FDR) in settings more realistic than i.i.d. (independent and identically
distributed), \citet{Kim08} use Gaussian deviates to simulate nested
correlation matrices with constrained correlation strengths.

In order to maximize expected performance, \citet{Nelson01} use a
stochastic simulation to compare a variety of management systems (e.g.,
queues). The authors use a modification of the method of \citet
{Marsaglia84} to simulate correlation matrices. Instead of generating
random points on the $k$-dimensional unit sphere, they constrain their
search to the part of the unit sphere with all nonnegative coordinates
(inducing nonnegative correlations).

In factor analysis, sample correlation matrices based on population
correlation matrices are typically used in simulation studies. Methods
which incorporate model error as well as sampling error create more
realistic structures from which to model data. \citet{Hong99}
recommends using the eigenstructure of the population correlation
matrix along with random chi-square deviates to directly obtain a
random sample covariance matrix, from which the sample correlation
matrix can be computed.

Additionally, as we detail in Section \ref{application}, simulated
correlation matrices are used to estimate heritability in Genome Wide
Association Studies (GWASs) and to assess network and classification
algorithms. We present work done by \citet{Lee11} on estimating
heritability without considering unknown noise in the correlation
structure. Alternatively, some clustering and classification methods
simulate correlations (or covariances) using uniform distributions
[\citet{Kraus10},
\citet{Tai07},
\citet{Huang10},
\citet{Kraj08}]. However, randomly simulated
univariate correlations do not typically produce matrices that are
positive semidefinite. We argue that neither the no-noise strategy nor
the uniform-noise strategy is realistic for assessing methodology. Our
work is applicable to any context where simulating realistic
correlation matrices is important.

Suppose that we are given a $N \times N$ correlation matrix $\Sigma=
(\Sigma_{ij})_{i,j=1}^N$. Generating a noisy correlation matrix $S =
(S_{ij})_{i,j=1}^N$ based upon the template $\Sigma$ can be difficult
since noise must be added to $\Sigma$ in such a way that $S$ remains
positive semidefinite and satisfies $S_{ii} =1$ and $-1 \leq S_{ij}
\leq1$ for $1 \leq i,j \leq N$. Moreover, for numerical purposes
(e.g., generating data from $S$) one might also require an explicit
upper bound on the condition number $\kappa(S)$ of $S$ (see Section
\ref
{SectionPreliminaries}) to ensure its numerical stability (e.g., for
matrix inversion). Unfortunately, naively adding random noise to a
correlation matrix can result in matrices which violate any or all of
the above constraints.

\subsection{Simulating data for evaluating algorithms}\label{simcor}

In certain applications, it is important to have a known (or assumed)
structure based on a covariance or correlation matrix. Many authors use
a particular structure and simulate Gaussian data from that matrix. For
example, in a recent paper, \citet{Tritchler09} simulate Gaussian data
to assess a method for filtering genes prior to, for example, network
analysis. Their structure consists of within group correlations of 0.4
and between group correlations of 0. Using clustering to find
differentially expressed genes, \citet{Hu10} generate Gaussian deviates
in a two-group setting with one group of 100 observations correlated at
0.94, another group of 608 observations clustered at 0.9, and
observations from different groups correlated at 0.697.

We appreciate the difficulty in generating realistic data with known
structure. However, we believe that using Gaussian deviates often adds
an additional unnecessary layer of assumptions.\vadjust{\goodbreak} Indeed, much recent
work has been applied to high-throughput data and, for example, we do
not believe that microarray data have Gaussian distributions [\citet
{Hardin2009}]. In Section \ref{discussion} we demonstrate that our
method produces matrices that are more general than the class of
matrices produced by finding the sample correlation of Gaussian data.
In particular, our method is able to produce Gaussian-like deviates or
other distributional deviates (e.g., uniform-like deviates). Instead of
simulating Gaussian data from a known correlation structure, we argue
in favor of simulating correlation matrices directly based on a known
correlation structure. The random correlation matrices can then be used
to assess the algorithm at hand.

\subsection{Three existing models}

The goal of our work is to provide an algorithm for simulating
correlation structures that can be used to evaluate statistical
methodology in a realistic context. Instead of relying on a known
structure, noise is added to the matrix to represent variability across
different components of the entries. Additionally, the noise is added
in a way that can represent any underlying data structure.

Below we have outlined three methods for generating correlation
matrices, each of which describes different dependence structures for
simulating data. Each of the three methods is taken from a different
area of application (estimating heritability from GWAS, classification
and network analysis). Our paper offers a flexible way to generate
correlation structures given a reasonable model of what we would expect
across observational units.

\subsubsection{Constant correlation model}
\label{const_corr_model}
Heritability is the proportion of variability in a phenotypic
characteristic which is due to genetic differences between individuals.
The estimation procedure for heritability is based on a mixture model
specified by a large correlation structure defining the correlations
between the genetic effects of individuals in a study. These
correlations are typically referred to as genetic correlations. The
genetic correlation structure is then used to decompose the phenotypic
variance to genetic and environmental components, resulting in an
estimate of heritability.

Recent work has assumed that the genetic correlation structure is known
[\citet{Lee11}], despite estimating it from genetic data. Simulations
using a known correlation structure are used to evaluate the
heritability estimation algorithm. In the actual data analysis, the
estimate of the correlation matrix is plugged into the algorithm as if
it were the true value.

The simulation study generates genetic relationships between 10,000
individuals in the following manner: Groups of size 100 are simulated
to have genetic correlations of 0.05. Uncorrelated environmental
effects are added to the genetic effects. The variances of the genetic
and environmental effects are predetermined by the value of
heritability used in the simulation. According to the liability
threshold model used by \citet{Lee11}, individuals for which the sum of
effects crosses a pre-defined threshold are considered to be cases,
while the rest are considered controls.

Since the prevalence of most interesting phenotypes is small, the
threshold is set such that only a small fraction of the individuals in
each group are considered cases. To simulate realistic case-control
studies, the cases are kept along with the same number of randomly
selected controls from the group, while the rest of the controls are
discarded. The process is repeated until 5000 cases and 5000 controls
are obtained.

Depending on the number of cases in each group, the resulting groups of
genetic correlation 0.05 are as small as a few individuals or as large
as 100 individuals. Genetic correlations between different groups are
assumed to be zero [\citet{Lee11}]. For future reference, we let
%
\begin{equation}
\label{eq-SigmaLee} \Sigma_k= \pmatrix{ 1 & 0.05 & 0.05 & 0.05 &
\cdots& 0.05
\cr
0.05 & 1 & 0.05 & 0.05 & \cdots& 0.05
\cr
0.05 & 0.05 & 1 & 0.05 &
\cdots& 0.05
\cr
0.05 & 0.05 & 0.05 & 1 & \cdots& 0.05
\cr
\vdots& \vdots& \vdots&
\vdots& \ddots& \vdots
\cr
0.05 & 0.05 & 0.05 & 0.05 & \cdots& 1}
\end{equation}
be the genetic correlation matrix for the $k$th group of individuals,
where the size of the $k$th block matrix is a random variable (i.e.,
the group size), with a distribution which is defined by the parameter settings.

\subsubsection{Toeplitz model} \label{subsubToep}

Another structure is one that models high correlation for observations
which are close together in the correlation matrix and models
decreasing correlation values for observations which are increasingly
far away. In building a classification model, \citet{guo07} describe a
Toeplitz structure (sometimes referred to as an auto-regressive
structure) to the correlation matrix, where adjacent pairs of
observations are highly correlated, and those further away are less
correlated. For future reference we let
%
\begin{equation}
\label{eq-SigmaToeplitz} \Sigma_k= \pmatrix{ 1 & \rho_k
& \rho_k^2 & \rho_k^3 & \cdots&
\rho_k^{g_k-1}
\cr
\rho_k & 1 &
\rho_k & \rho_k^2 & \cdots&
\rho_k^{g_k-2}
\cr
\rho_k^2 &
\rho_k & 1 & \rho_k & \cdots& \rho_k^{g_k-3}
\cr
\rho_k^3 & \rho_k^2 &
\rho_k & 1 & \cdots& \rho_k^{g_k-4}
\cr
\vdots&
\vdots& \vdots& \vdots& \ddots& \vdots
\cr
\rho_k^{g_k-1} &
\rho_k^{g_k-2} & \rho_k^{g_k-3} &
\rho_k^{g_k-4} & \cdots& 1}
\end{equation}
be the correlation matrix for the $k$th class, given by the base
correlation value $\rho_k$. In this model, the between group
correlations are set to zero. Additional classification models have
used similar Toeplitz structure for simulating data from a correlation
matrix [\citet{guo07},
\citet{Storey07},
\citet{Witten09},
\citet{Zuber09},
\citet{Pang09},
\citet{Huang10}]. In
fact, Huang et al. use a $U[0.5,1.5]$ distribution to simulate the
variance components in order to add noise to the above prescribed
structure. Further, this Toeplitz correlation structure is seen in time
series models where simulating correlation matrices is also important
[\citet{Joe06},
\citet{Ng10}].
The Toeplitz structure has been used extensively in classification and
discriminant analysis as a model for group correlations

\subsubsection{Hub observation model} \label{subsubhub}

The last model which we consider is one that is hierarchical in nature
based on a single hub-observation and the relationship of each
observation to that original hub. Within the context of network
analysis, Horvath et al. [\citet{horvath05},
\citet{horvath08a},
\citet{horvath08b}]
define a structure with respect to a particular profile (or
hub-observation). Each observation in a group is correlated with the
hub-observation with decreasing strength (from a user supplied maximum
correlation to a given minimum correlation). Additionally, groups are
generated independently (i.e., with correlation zero between
groups). Letting observation 1 correspond to the hub, for the $i$th
observation ($i=2, 3, \ldots, g$), the correlation between it and the
hub-observation is given by
\[
\Sigma_{i,1} = \rho_{\max} - \biggl(\frac{i-2}{g-2}
\biggr)^{\gamma}(\rho _{\max} - \rho_{\min}).
\]
Note that the correlation between the $i$th observation and the hub
will range from $\rho_{\max}$ to $\rho_{\min}$; the rate at which the
correlations decay is controlled by the exponent~$\gamma$ (where
$\gamma= 1$ would indicate a linear decay).

\subsubsection{Overview}
Motivated by the models above, we provide algorithms for adding noise
to prescribed correlation matrices. We begin in Section~\ref
{SectionExamples} detailing algorithms for the three specific models of
correlation matrices discussed above. In Section \ref{discussion} we
demonstrate the benefits of generating random deviates from the
correlation matrix instead of using random deviates from a particular
distribution. Section \ref{application} gives applications of how our
method can be used to assess new and standard statistical procedures.
Following a brief conclusion in Section \ref{SectionConclusion}, we
present the theoretical justifications of our algorithms in the \hyperref[Appendix]{Appendix}.

\section{Recipes}\label{SectionExamples}

Using a single basic procedure (Algorithm \ref{AlgorithmMain} in
Section \ref{SubsectionBasic}) for adding noise to a given
correlation matrix, we can take advantage of our theoretical
understanding of certain known correlation structures to yield stronger
results. This is carried out for the constant correlation structure
(Algorithm \ref{AlgorithmBlocks} in Section~\ref{ConstantCorrelationAlg}),
the Toeplitz correlation structure (Algorithm~\ref{AlgorithmToeplitz}
in Section~\ref{SubsectionToeplitzAlg}) and the hub correlation
structure (Algorithm \ref{AlgorithmToeplitzHub} in Section~\ref
{SubsectionHubAlg}).
Each model describes a population based on multiple groups with the
same underlying structure (with different sizes and parameter values).
Since the justifications of these procedures are rather involved, we
defer the technical details until \hyperref[Appendix]{Appendix}.

\subsection{Constant correlation structure}\label{ConstantCorrelationAlg}

Our first correlation structure is\break based on constant correlations
within each group and between each group (values of the correlation
differ for each relationship). In particular, observe that the approach
below yields a noisy correlation matrix which has a significant amount
of noise on the \emph{off-diagonal} blocks. This is clearly more
realistic than simply assuming that all of these entries are zero.
A detailed justification of the following algorithm can be found in
Section \ref{SubsectionAlgorithmBlocks} of
\hyperref[Appendix]{Appendix}.

\begin{Algorithm}\label{AlgorithmBlocks}
Let
\begin{itemize}
\item$K$ denote a positive integer (the number of groups) and $k =
1,2,\ldots, K$,
\item$g_k$ be a positive integer (the size of the $k$th group),\vspace*{1pt}
\item$N = \sum_{k=1}^K g_k$ (size of the desired matrix),
\item$\rho_k$ such that $0 \leq\rho_k < 1$ (baseline correlation in
the $k$th group),
\item$\rho_{\min} = \min\{\rho_1,\rho_2,\ldots,\rho_K\}$ (minimum
correlation in any group),
\item$\rho_{\max} = \max\{\rho_1,\rho_2,\ldots,\rho_K\}$ (maximum
correlation in any group),
\item$\delta$ such that $0 \leq\delta< \rho_{\min}$ (baseline noise
between group),
\item$\Sigma_k$ be the $g_k \times g_k$ matrix
%
\begin{equation}
\label{eq-SigmaBlock} \Sigma_k= \pmatrix{ 1 & \rho_k &
\cdots& \rho_k
\cr
\rho_k & 1 & \cdots&
\rho_k
\cr
\vdots& \vdots& \ddots& \vdots
\cr
\rho_k &
\rho_k & \cdots& 1}
\end{equation}
(correlation matrix for $k$th group),
\item$\Sigma$ be the $N \times N$ matrix having the blocks
$\Sigma_1, \Sigma_2, \ldots, \Sigma_k$ along the diagonal and zeros
elsewhere,
\item$\epsilon$ such that $0 \leq\epsilon< 1 - \rho_{\max}$ (maximum
entry-wise random noise),
\item$M$ be a positive integer (the dimension of the noise space).
\end{itemize}
Select $N$ unit vectors $\vec{u}_1, \vec{u}_2, \ldots, \vec{u}_N$
randomly from $\R^M$.
The $N \times N$ matrix $S = (S_{ij})_{i,j=1}^N$ defined by
%
\begin{equation}
\label{eq-MatrixBlocks} S_{ij} = \cases{ 1, &\quad $\mbox{if $i =j$},$
\vspace*{2pt}\cr
\rho_k + \epsilon\vec{u}_i^T
\vec{u}_j, & \quad $\mbox{if $i,j$ are in the $k$th group and $i \neq j$},$
\vspace*{2pt}\cr
\delta+ \epsilon\vec{u}_i^T \vec{u}_j,&\quad
$\mbox{if $i,j$ are in different groups},$}
\end{equation}
is a correlation matrix whose condition number satisfies
%
\begin{equation}
\label{eq-KappaBlocks} %
\fbox{
$\displaystyle\kappa(S) \leq\frac{N(1+\epsilon)+1}{1-\rho_{\max}-\epsilon}.$
}
\end{equation}
\end{Algorithm}

\subsection{Toeplitz correlation structure}\label{SubsectionToeplitzAlg}

The Toeplitz structure has been used extensively in classification,
discriminant analysis and in the time series literature as a model for
group correlations. In particular, the model we follow assumes that
each pair of adjacent observations is highly correlated and that the
correlations between the $i$th and $j$th observations decay
exponentially with respect to $|i-j|$. The following algorithm, whose
justification can be found in Section \ref
{SubsectionAlgorithmToeplitz} of \hyperref[Appendix]{Appendix}, produces
noisy correlation matrices based upon the Toeplitz template.

\begin{Algorithm}\label{AlgorithmToeplitz}
Let
\begin{itemize}
\item$K$ denote a positive integer (the number of clusters) and $k =
1,2,\ldots, K$,
\item$g_k$ be a positive integer (the size of the $k$th group),
\item$N = \sum_{k=1}^K g_k$ (size of the desired matrix),
\item$\rho_k$ be such that $0 \leq\rho_k < 1$ (correlation factor in
the $k$th group),
\item$\rho_{\max} = \max\{\rho_1,\rho_2,\ldots,\rho_K\}$ (maximum
correlation factor),
\item$\Sigma_k$ be the $g_k \times g_k$ Toeplitz correlation matrix
%
\begin{equation}
\Sigma_k= \pmatrix{ 1 & \rho_k &
\rho_k^2 & \rho_k^3 & \cdots&
\rho_k^{g_k-1}
\cr
\rho_k & 1 & \rho_k
& \rho_k^2 & \cdots& \rho_k^{g_k-2}
\cr
\rho_k^2 & \rho_k & 1 &
\rho_k & \cdots& \rho_k^{g_k-3}
\cr
\rho_k^3 & \rho_k^2 &
\rho_k & 1 & \cdots& \rho_k^{g_k-4}
\cr
\vdots&
\vdots& \vdots& \vdots& \ddots& \vdots
\cr
\rho_k^{g_k-1} &
\rho_k^{g_k-2} & \rho_k^{g_k-3} &
\rho_k^{g_k-4} & \cdots& 1}
\end{equation}
(correlation matrix for $k$th group),
\item$\Sigma$ be the $N \times N$ matrix having the blocks
$\Sigma_1, \Sigma_2, \ldots, \Sigma_k$ along the diagonal and zeros
elsewhere,

\item$0 < \epsilon< \dfrac{1 - \rho_{\max}}{1 + \rho_{\max}}$
(maximum entry-wise random noise),

\item$M$ be a positive integer (the dimension of the noise space).
\end{itemize}
Select $N$ unit vectors $\vec{u}_1, \vec{u}_2, \ldots, \vec{u}_N$ from
$\R^M$
and form the $M \times N$ matrix $U=(\vec{u}_1 | \vec{u}_2 | \cdots|
\vec{u}_N)$ whose columns
are the $\vec{u}_i$. The $N \times N$ matrix
%
\begin{equation}
\label{eq-MatrixToeplitz} S = \Sigma+\epsilon\bigl(U^TU - I\bigr)
\end{equation}
is a correlation matrix whose entries satisfy $|S_{ij} - \Sigma_{ij} |
\leq\epsilon$
and whose condition number satisfies
%
\begin{equation}
\label{eq-KappaToeplitz} %
\fbox{
$\displaystyle\kappa(S) \leq \frac{ {(1+ \rho_{\max})}/{(1-\rho_{\max})} + (N-1)\epsilon}{
{(1-\rho_{\max})}/{(1+\rho_{\max})}-\epsilon}$.}
\end{equation}
\end{Algorithm}

Among other things, let us remark that for typical values of $\rho$
[e.g., \citet{guo07} let $\rho= 0.9$] the noise level $\epsilon$ can
be made quite large compared to most of the entries in each $\Sigma_k$.
This occurs because the eigenvalue estimates \eqref{eq-ToeplitzSharp}
obtained in Section \ref{SubsectionAlgorithmToeplitz} are remarkably
strong and because the off-diagonal entries of each submatrix $\Sigma
_k$ are small (due to exponential decay) if one is far away from the
main diagonal. Thus, the approach outlined above yields a flexible
method for introducing noise into the Toeplitz model. In fact, one can
introduce so much noise (while still obtaining a correlation matrix
with controlled condition number) that the original block-Toeplitz
structure becomes difficult to discern.

\subsection{Hub correlation structure}\label{SubsectionHubAlg}

The hub correlation structure assumes a known correlation between a
\emph{hub observation} (typically the first observation) and each of
the other observations. Moreover, one typically assumes that the
correlation between the $1$st and the $i$th observation decays as $i$ increases.

Let us describe a typical example which has been considered frequently
in the literature.
Suppose that the first row (and hence column) of a $g \times g$
correlation matrix $A$ is to
consist of the prescribed values
\[
A_{11} = 1,\qquad A_{1i} = \rho_{\max} - (
\rho_{\max} - \rho_{\min}) \biggl( \frac
{i-2}{g-2}
\biggr)^{\gamma}
\]
which decrease (linearly if $\gamma=1$) from $A_{12} = \rho_{\max}$ to
$A_{1g} = \rho_{\min}$ for $2 \leq i \leq g$.
For instance, this model is considered in Horvath et al.~[\citet
{horvath05,horvath08a,horvath08b}]. For the sake of simplicity,
we consider the linear case $\gamma=1$ and adopt a more convenient
notation. Rather than specifying $\rho_{\max}$
and $\rho_{\min}$, we specify only $\rho_{\max}$ and work instead with
the step size $\tau= (\rho_{\max} - \rho_{\min})/(g-2)$.

After specifying the first row, there are a variety of ways to generate
the remainder of such a correlation matrix. Using any hub structure
correlation matrix, we can find the smallest resulting eigenvalue which
can be fed into Algorithm \ref{AlgorithmMain} of Section~\ref
{SubsectionBasic}. For example, we can use a Toeplitz structure to fill
out the remainder of the hub correlation matrix and, using the
well-developed theory of truncated Toeplitz matrices [\citet{LTTM}],
obtain eigenvalue bounds which can be fed directly into Algorithm \ref
{AlgorithmMain}. This approach yields the following algorithm, whose
justification can be found in Section \ref
{SubsectionAlgorithmToeplitzHub} of \hyperref[Appendix]{Appendix}.

\begin{Algorithm}\label{AlgorithmToeplitzHub}
Let
\begin{itemize}
\item$K$ denote a positive integer (the number of groups) and $k =
1,2,\ldots, K$,
\item$g_k$ be a positive integer (the size of the $k$th group),
\item$N = \sum_{k=1}^K g_k$ (size of the desired matrix),
\item$\rho_k$ (maximum correlation in the first row of $k$th group),
\item$\tau_k$ (step size in first row/column of $k$th group),
\item$\alpha_{k,1} = 1$ and $\alpha_{k,i} = \rho_k - \tau_k (i-2)$
(correlations between hub and observations),\eject
\item$\Sigma_k$ be the $g_k \times g_k$ hub-Toeplitz correlation matrix
%
\begin{equation}
\label{eq-SigmaToeplitzHub} \Sigma_k= \pmatrix{ 1 &
\alpha_{k,2} & \alpha_{k,3} & \alpha_{k,4} & \cdots&
\alpha _{k,g_k}
\cr
\alpha_{k,2} & 1 & \alpha_{k,2} &
\alpha_{k,3} & \cdots& \alpha _{k,g_k-1}
\cr
\alpha_{k,3} &
\alpha_{k,2} & 1 & \alpha_{k,2} & \cdots& \alpha _{k,g_k-2}
\cr
\alpha_{k,4} & \alpha_{k,3} & \alpha_{k,2} & 1 &
\cdots& \alpha _{k,g_k-3}
\cr
\vdots& \vdots& \vdots& \vdots& \ddots& \vdots
\cr
\alpha_{k,g_k} & \alpha_{k,g_k-1} & \alpha_{k,g_k-2} & \alpha
_{k,g_k-3} & \cdots& 1}
\end{equation}
(correlation matrix for $k$th group),
\item$\Sigma$ be the $N \times N$ matrix having the blocks
$\Sigma_1, \Sigma_2, \ldots, \Sigma_k$ along the diagonal and zeros
elsewhere,

\item$0 < \epsilon< \min \{1 - \rho_k - \frac{3}{4} \tau_k \dvtx 1
\leq k \leq K  \}$ ($\epsilon$ is the maximum noise level),

\item$M$ be a positive integer (the dimension of the noise space).
\end{itemize}
Select $N$ unit vectors $\vec{u}_1, \vec{u}_2, \ldots, \vec{u}_N$ from
$\R^M$
and form the $M \times N$ matrix $U=(\vec{u}_1 | \vec{u}_2 | \cdots|
\vec{u}_N)$ whose columns
are the $\vec{u}_i$. The $N \times N$ matrix
%
\begin{equation}
\label{eq-MatrixToeplitzHub} S = \Sigma+\epsilon\bigl(U^TU - I\bigr)
\end{equation}
is a correlation matrix whose entries satisfy $|S_{ij} - \Sigma_{ij} |
\leq\epsilon$
and whose condition number satisfies
%
\begin{equation}
\label{eq-BasicKappa1} %
\fbox{
$\displaystyle\kappa(S) \leq\frac{ \lambda_1(\Sigma) + (N-1)\epsilon}{ \lambda
_N(\Sigma) - \epsilon} $}
\end{equation}
where
%
\begin{eqnarray}
\lambda_1(\Sigma) &\leq&\max \biggl\{1 + (g_k-1)
\rho_k - \tau_k \frac
{(g_k-2)(g_k-1)}{2}\dvtx 1\leq k \leq K
\biggr\}, \label{eq-EVBig}
\\
\lambda_N(\Sigma) &\geq&\min \bigl\{1 - \rho_k -
\tfrac{3}{4}\tau _k \dvtx 1 \leq k \leq K \bigr\}.
\label{eq-EVSmall}
\end{eqnarray}
\end{Algorithm}

\subsection{Extensions}

Before proceeding, let us remark that our general Algorithm~\ref
{AlgorithmMain},
which can be found in Section \ref{SubsectionBasic} of \hyperref[Appendix]{Appendix},
is applicable to any given positive-definite correlation
matrix. The amount of noise which can be added to the original matrix
is determined by its smallest eigenvalue. For several specific classes
of correlation matrices, one can obtain simple, but powerful, lower
bounds on this lowest eigenvalue. For such correlation matrices, we
have provided explicit, specialized algorithms which provide a
significant amount of noise while also maintaining quantitative control
over the condition number of the resulting matrix.

\section{Distribution of error terms}\label{disterror}

As described above, our method uses the dot product of normalized
vectors as the error terms which are added to a given correlation
matrix. Below we discuss three methods for generating normalized
vectors with given distributions.

\begin{longlist}[1.]
\item[1.]
\textit{Random uniform vectors on the M-dimensional unit sphere}: Consider
\begin{eqnarray*}
x_i &\sim& \mbox{i.i.d. } N(0,1),  \qquad  i=1, 2, \ldots, M,
\\
\mathbf{x} &=& (x_1, x_2, \ldots, x_M),
\\
\mathbf{v} &=& \frac{\mathbf{x}}{ \Vert \mathbf{x}\Vert  }.
\end{eqnarray*}
It is known that $\mathbf{v}$ will be uniformly distributed on the
$M$-dimensional unit sphere [\citet{Muller59}]. Additionally, for
vectors distributed uniformly on the unit sphere, the distribution of
their dot product is well characterized [\citet{Cho09}]:
\begin{eqnarray*}
\mathbf{v}, \mathbf{w} &\sim& \mbox{uniformly on the $M$-dimensional unit
sphere},
\\
Z &=& \mathbf{v}^T \mathbf{w},
\\
f_Z(z;M) &=& \frac{\Gamma({M}/{2})}{\Gamma({(M-1)}/{2})
\sqrt{\pi
}} \bigl(\sqrt{1-z^2}
\bigr)^{M-3},\qquad   -1 \leq z \leq1,
\end{eqnarray*}
is the probability density function for the dot product of $\mathbf{v}$ and
$\mathbf{w}$. Note that if $M=2$, the distribution of $Z$ is of the form $1/
\pi\sqrt{1-z^2}$ which gives a U-shaped distribution favoring values
of $Z$ closer to $-1$ and 1. If $M=3$, the distribution is uniform
across $-1$ to 1. For $M > 3$, the distribution function is mound
shaped and converges to a Gaussian distribution for large $M$ (see below).

\item[2.]
\textit{Random independent and identically distributed vectors}: Consider
two vectors generated independently from identical distributions in $\R
^M$ with mean zero,
\begin{eqnarray*}
\mathbf{v}, \mathbf{w} &\sim& \mbox{i.i.d. } F_M ({\bf\mu}=
\mathbf{0}),
\\
Z &=& \frac{\mathbf{v}^T \mathbf{w}}{\Vert \mathbf{v}\Vert  \Vert \mathbf{w}\Vert },
\\
\sqrt{M} Z &\stackrel{D} {\rightarrow}& N(0,1).
\end{eqnarray*}
The asymptotic distribution of $Z$ is a straightforward application of
the Central Limit theorem and Slutsky's theorem.

\item[3.]
\textit{Arbitrary distribution}: Some situations may call for a particular
distribution of error noise. The distribution can be controlled through
the $\alpha$ parameter as seen in equation (\ref{eq-SquareRoot}).
\end{longlist}

Note that typically the error terms added to the correlation entries
are of the form
\[
\mbox{error} = \epsilon\cdot\frac{\mathbf{v}^T \mathbf{w}}{\Vert \mathbf{v}\Vert
\Vert \mathbf{w}\Vert }.
\]

If the dot product is approximately distributed with a variance of
$1/M$, then the variance of the error term is $\epsilon^2/M$, resulting
in a standard error of the correlation values,
%
\begin{equation}
\label{SEcor} \operatorname{SE}( \mbox{correlation} ) \approx\frac{\epsilon}{\sqrt{M}}.
\end{equation}

The distribution of error terms will necessarily depend on the
application. For some problems, uniform error terms may be most
appropriate; for other problems, Gaussian errors will be preferable. In
fact, for Gaussian data, correlations between vectors are approximately
Gaussian, which may motivate a user to want to add Gaussian noise to
the given correlation structure.

\subsection{Comparison to a correlation matrix from Gaussian
vectors}\label{discussion}

One method for generating a noisy correlation matrix is to simulate
Gaussian data from an original template and then find the sample
correlation matrix from the data. Varying the sample size of the
generated data can create correlation matrices which are more or less
variable (in magnitude). However, from Gaussian data the nature of the
variability (distribution) of the resulting correlations is similar
across different sample sizes---uniform or U-shaped distributions of
error terms are not possible given correlations from Gaussian data. In
addition, the majority of the entries in a given sample correlation
matrix generated from Gaussian data are quite close to the template
matrix. Only a handful of observations deviate from the template
substantially. In fact, the sample size needed in order to get a large
amount of variability could be smaller than the dimension of the
correlation matrix (thus producing sample correlation matrices which
are not positive definite).

\begin{figure}

\includegraphics{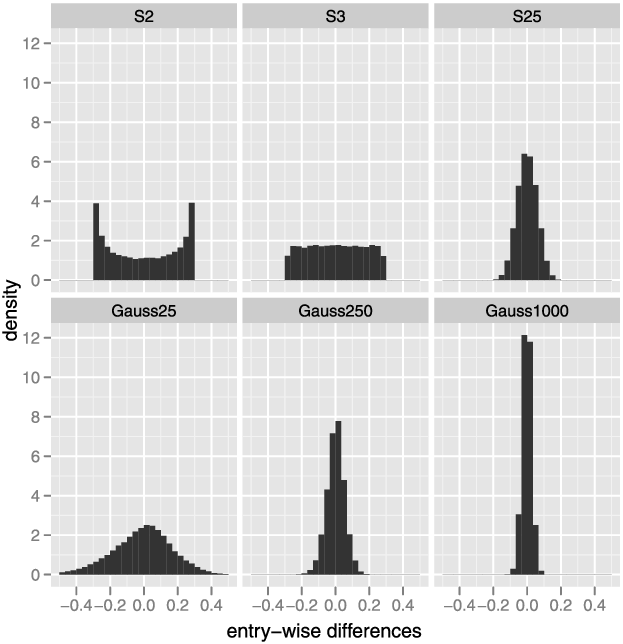}

\caption{Each histogram represents the distribution of entry-wise
differences between the generated matrix and the template. The
distribution of differences for random
vectors in $\R^{25}$ (S25) is similar to that obtained from the
correlation matrix generated
by sampling 250 random Gaussian vectors (Gauss250).}\label{histNormsim}\vspace*{-2pt}
\end{figure}

To demonstrate the restriction associated with simulating Gaussian data
as a way to find sample correlation matrices, we generate multiple
correlation matrices using both Gaussian samples and our method. The
Gaussian noise is created by simulating data of a particular sample
size (25, 250, or 1000) from a template correlation matrix. We then
compute the sample correlation matrix and find the difference between
the estimate and the template; histograms of those differences describe
the distribution of the correlation error terms. For example, Gauss25
was created by simulating 25 observations from a $230 \times230$
template correlation matrix. The difference between the correlations of
the 25 observations and the template matrix are computed; the
histograms of the differences are given in Figure~\ref{histNormsim}.

The three Gaussian structures (sample sizes 25, 250, 1000) show the
same tendencies with more spread for smaller sample sizes (see Figure
\ref{histNormsim}: Gauss25, Gauss250, Gauss1000). The three simulations
using our method are based on uniform random vectors on the unit sphere
in $\R^{2}$, $\R^{3}$ and $\R^{25}$ (see Table~\ref{NormSim} for different generating scenarios).
For each simulation we used a
constant correlation structure with three groups of sizes $g_1=100,
g_2=50, g_3=80$; within group correlations of $\rho_1=0.7, \rho_2=0.7,
\rho_3=0.4$; $\epsilon=0.29$; and between group correlations of
$\delta
=0.25$ (see Algorithm \ref{AlgorithmBlocks} in Section \ref
{SubsectionAlgorithmBlocks}).

Note that in our simulations below, the distribution of differences
from uniform vectors in $\R^{25}$ (S25) is similar to that obtained
from the correlation matrix generated by sampling 250 random Gaussian
vectors (Gauss250). In fact, not only can our method create
correlations with univariate distributions that are similar to Gaussian
deviations, but the eigenvalues of the respective matrices (ours
compared with correlations from Gaussian data) are indistinguishable
(results not shown).

We see that our method is able to add larger noise terms than the
Gaussian simulation. Figure \ref{histNormsim} shows the distribution of
the differences. Depending on the application, one might prefer large
noise components, uniform noise components or bell-shaped noise
components. Our work provides a template for generating multiple
different structures depending on the problem at hand.

\begin{table}
\caption{Six different correlation matrix
generating scenarios. \textup{S2}, \textup{S3} and \textup{S25} use the algorithms given in the
paper for constant correlation; to simulate the error terms, we
generate vectors on the unit sphere using standard Gaussian deviates.
The Gaussian simulations use the template matrix with the given sample
size of random vectors. Each correlation matrix is based on a setting
of 3 groups with sizes $(100, 50, 80)$}\label{NormSim}
\begin{tabular*}{\textwidth}{@{\extracolsep{\fill}}lccccc@{}}
\hline
{Template} &
$\rho_1= 0.7$ & $\rho_2 = 0.7$ & $\rho_3=0.4$ &
$\delta= 0.25$ & $\epsilon= 0$
\\[3pt]
\mbox{S2} & $\rho_1= 0.7$ & $\rho_2 = 0.7$ &
$\rho_3=0.4$ & $\delta= 0.25$ & $\mathbf{u}_i \in
S^2, \epsilon= 0.29$
\\
\mbox{S3} & $\rho_1= 0.7$ & $\rho_2 = 0.7$ &
$\rho_3=0.4$ & $\delta= 0.25$ & $\mathbf{u}_i \in
S^3, \epsilon= 0.29$
\\
\mbox{S25} & $\rho_1= 0.7$ & $\rho_2 = 0.$7 &
$\rho_3=0.4$ & $\delta= 0.25$ & $\mathbf{u}_i \in
S^{25}, \epsilon= 0.29$
\\[3pt]
\mbox{Gauss25} & $\rho_1= 0.7$ & $\rho_2 = 0.7$ &
$\rho_3=0.4$ & $\delta= 0.25$ & \mbox{25 vectors}
\\
\mbox{Gauss250} & $\rho_1= 0.7$ & $\rho_2 = 0.7$ &
$\rho_3=0.4$ & $\delta= 0.25$ & \mbox{250 vectors}
\\
\mbox{Gauss1000} & $\rho_1= 0.7$ & $\rho_2 = 0.7$&
$\rho_3=0.4$ & $\delta= 0.25$ & \mbox{1000 vectors}
\\
\hline
\end{tabular*}
\end{table}

\section{Applications} \label{application}

To demonstrate the effectiveness of our method, we simulate data from
two applications to show that noise added to a known correlation
structure can be useful in practice. It may not always be obvious which
format to use to incorporate the noise; the format of the noise will be
situation dependent and should be based on the underlying data
structure. In Section \ref{discussion} we have provided more details
about the different noise models.

\subsection{Heritability}\label{SubsectionHerit}

Heritability is the proportion of variability in a phenotypic
characteristic which is due to genetic differences. The understanding
and estimation of heritability is of great importance in directing
future studies as well as understanding the architecture of human
genetic diseases such as type-1 and type-2 diabetes, Crohn's disease,
schizophrenia and bipolar disorder. The study of heritability in human
disease presents the so-called ``mystery'' of the missing heritability
[\citet{Maher08}]: a considerable gap between the estimated heritability
as obtained from family studies and the estimated heritability as
obtained from genetic studies (known as genome-wide association
studies, or GWASs), with the latter estimate of heritability being
considerably smaller than the former for a wide range of phenotypes.

Recently, a novel method for estimation of heritability from the
genotypes of unrelated individuals was introduced by \citet{Yang10}.
The method was first applied for random samples from a population
[\citet
{Yang10}] and later adapted for the more relevant scenario of
case-control studies [\citet{Lee11}]. These works presented genome-based
estimates of heritability that were considerably higher than previous
estimates, thus bridging, at least in part, the gap between
family-based and genome-based estimates of heritability. As expected,
these works attracted attention and are the focus of recent research
and debate [see, e.g., \citet{golan11},
\citet{Lee12}].

The central idea behind these methods is to estimate a population-wise
correlation structure from the genotypes of individuals and use this
estimated structure in a Restricted Maximum Likelihood (REML)
estimation of the heritability. However, the REML estimation does not
account for the fact that the correlation structure is estimated rather
than known. Moreover, the simulations in \citet{Lee11} use a known
correlation structure to demonstrate the validity of the method, which
in turn uses an estimated correlation structure. Such simulations might
produce an overly optimistic evaluation of the method used to estimate
heritability in terms of both bias and variance.

To demonstrate the sensitivity of the heritability estimate to the
known correlation structure, we reran the simulations in \citet{Lee11}
with and without noise. As expected, adding noise to the matrix
introduces bias to the estimators. Our methods provide a mechanism for
understanding the behavior of heritability estimates under different
correlation and error structures.

For our investigation, we are interested in estimating heritability in
the setting of a binary response, in particular, we want to estimate
heritability for case control studies. As done in \citet{Lee11}, we
assume there is an underlying liability continuous variable (e.g.,
glucose level) determining the binary measured phenotype (e.g.,
diabetes). We can find the heritability on the observed scale (with
respect to the binary disease trait) and transform it back to the value
of interest, the heritability on the liability scale. The
transformation considers the disease prevalence in the population
[\citet
{Lee11}].

We followed the simulation procedure of \citet{Lee11}, which is
outlined in Section~\ref{const_corr_model}. See \citet{Lee11} for a
more detailed description of the simulation procedure. Subsequently, we
added noise to the correlation matrices using Algorithm \ref
{AlgorithmBlocks} in Section \ref{ConstantCorrelationAlg}. We used the
software Genome-wide Complex Trait Analysis (GCTA) to estimate
heritability and standard errors of the estimate [\citet{Yang11}].

Results are presented in Table \ref{TableBias}. Each table entry
contains the estimated heritability of liability from our simulations
with noise as well as the corresponding estimate given by \citet{Lee11}---calculated
using the known correlation structure. As expected, the
more noise added to the relationship matrix, the more bias in
estimating the heritability. Additionally, we see that there is a
strong interaction: for low prevalence, even a small amount of error
can have a large impact on the estimate of heritability. With high
prevalence, moderate amounts of error can bias the estimate.

\begin{table}
\tabcolsep=0pt
\caption{The average heritability from 100 simulations for the given
population prevalence and heritability of liability.
In parentheses we provide the estimate given by \citet{Lee11} with no
error. The error is added as described in Algorithm \protect\ref
{AlgorithmBlocks}
with $\epsilon= 0.001, 0.01$ and 0.02, corresponding to a SE of the
noise terms of 0.0002, 0.002 and 0.004, respectively; see equation
\protect\eqref{SEcor}.}\label{TableBias}
\begin{tabular*}{\textwidth}{@{\extracolsep{\fill}}lcd{1.3}ccccc@{}}
\hline
& & &\multicolumn{5}{c@{}}{\textbf{True heritability of liability}}\\[-6pt]
& & &\multicolumn{5}{c@{}}{\hrulefill}\\
&& \multicolumn{1}{c}{\textbf{Prevalence of disease in pop.}} & \textbf{0.1} & \textbf{0.3} & \textbf{0.5} & \textbf{0.7} & \textbf{0.9}\\
\hline
\multirow{3}{*}{\rotatebox{-90}{\fontsize{8.3}{10.3}{\selectfont{$\epsilon= 0.001$}}}} &
\multirow{3}{*}{\rotatebox{-90}{\fontsize{8.3}{10.3}{\selectfont{$\operatorname{SE} =
0.0002$}}}}
& 0.5 & 0.10 (0.09) & 0.28 (0.28) & 0.47 (0.51) & 0.64 (0.70) & 0.80
(0.90)\\[3pt]
&& 0.1 & 0.11 (0.11) & 0.29 (0.30) & 0.49 (0.49) & 0.70 (0.71) & 0.87
(0.89)\\[3pt]
&& 0.001 & 0.05 (0.17) & 0.23 (0.31) & 0.37 (0.56) & 0.56 (0.75) & 0.77
(0.94)\\[8pt]
\multirow{3}{*}{\rotatebox{-90}{\fontsize{8.3}{10.3}{\selectfont{$\epsilon= 0.01$}}}} & \multirow
{3}{*}{\fontsize{8.3}{10.3}{\selectfont{\rotatebox{-90}{$\operatorname{SE} = 0.002$}}}}
& 0.5 & 0.08 (0.09) & 0.23 (0.28) & 0.37 (0.51) & 0.54 (0.70) & 0.69
(0.90)\\[3pt]
&& 0.1 & 0.06 (0.11) & 0.23 (0.30) & 0.42 (0.49) & 0.60 (0.71) & 0.78
(0.89)\\[3pt]
&& 0.001 & 0.00 (0.17) & 0.01 (0.31) & 0.01 (0.56) & 0.02 (0.75) & 0.02
(0.94)\\[8pt]
\multirow{3}{*}{\rotatebox{-90}{\fontsize{8.3}{10.3}{\selectfont{$\epsilon= 0.02$}}}} & \multirow
{3}{*}{\fontsize{8.3}{10.3}{\selectfont{\rotatebox{-90}{$\operatorname{SE} = 0.004$}}}}& 0.5 & 0.07 (0.09) & 0.19 (0.28) & 0.35 (0.51) & 0.51 (0.70) & 0.65
(0.90)\\[3pt]
&& 0.1 & 0.03 (0.11) & 0.16 (0.30) & 0.22 (0.49) & 0.53 (0.71) & 0.73
(0.89)\\[3pt]
&& 0.001 & 0.01 (0.17) & 0.01 (0.31) & 0.01 (0.56) & 0.01 (0.75) & 0.01
(0.94)\\
\hline
\end{tabular*}
\end{table}

We run a second set of simulations to see how often we can capture the
true heritability in a confidence interval using $\pm$ 2 SE (provided
from the GCTA software) when noisy correlation matrices are used; we
expect roughly 95\% of the confidence intervals to capture the true
heritability value. We simulated 100 heritability values from a uniform $(0.1, 0.9)$ distribution. We then simulated phenotypes and a
corresponding correlation matrix with error (as described previously)
for each heritability. We constructed corrected CIs, using a
multiplicative factor correction obtained from our first set of
simulations as well as uncorrected confidence intervals, and counted
the number of times these CIs contained the true heritability.

Results for the CIs are given in Table \ref{TableCI}. Each entry gives
the number of true heritabilities captured in the interval two standard
errors around the adjusted estimated heritability (expected to be 95
when CIs are accurate). Our results show that neglecting to account for
the effects of noise on heritability estimation resulted in problematic
confidence intervals. Correcting the bias using the method described
above resolved the issue for low and moderate levels of noise.

\begin{table}
\caption{Out of 100 simulations, the number of true
heritabilities captured in the interval two standard errors around the
estimated heritability. We would expect 95\% confidence intervals to
capture the true heritability value 95 times out of 100. For each
parameter setting, we calculated two intervals, thus, we report two
coverage rates. The first number uses a correction factor (for both the
heritability and the SE of the heritability) calculated from the bias
estimated from Table \protect\ref{TableBias} above. The second number
uses no
correction. The error is added as described in Algorithm \protect\ref
{AlgorithmBlocks} with $\epsilon= 0.001, 0.01$ and 0.02, corresponding
to a SE of the noise terms of 0.0002, 0.002 and 0.004, respectively;
see equation (\protect\ref{SEcor}).}\label{TableCI}
\begin{tabular*}{\textwidth}{@{\extracolsep{\fill}}lcccc@{}}
\hline
& \multicolumn{4}{c@{}}{\textbf{Amount of error added to the
correlation structure}}\\[-6pt]
& \multicolumn{4}{c@{}}{\hrulefill}\\
 & \multicolumn{1}{c}{$\bolds{\epsilon= 0}$} & \multicolumn{1}{c}{$\bolds{\epsilon= 0.001}$} &
 \multicolumn{1}{c}{$\bolds{\epsilon= 0.01}$}&
\multicolumn{1}{c@{}}{$\bolds{\epsilon= 0.02}$}\\
\multicolumn{1}{@{}l}{\multirow{2}{90pt}[11pt]{\textbf{Prevalence of disease in population}}} & \multicolumn{1}{c}{$\bolds{\operatorname{SE} = 0}$} &
\multicolumn{1}{c}{$\bolds{\operatorname{SE} = 0.0002}$} &
 \multicolumn{1}{c}{$\bolds{\operatorname{SE} = 0.002}$}&
\multicolumn{1}{c@{}}{$\bolds{\operatorname{SE}=0.004}$}\\
\hline
0.5 & $93/89$ & $94/88$ & $94/76$ & $93/46$ \\
0.1 & $89/90$ & $94/94$ & $91/72$ & $68/57$ \\
\hline
\end{tabular*}
\end{table}

We note that for the case of estimating heritability, the bias suggests
an immediate way of choosing $\epsilon$ for proper inference. The bias
is a result of the noisy estimation of the true underlying genetic
correlation structure, and so can be estimated using the actual
genotypes, as done by \citet{Yang10}.
We thus suggest choosing $\epsilon$ such that the bias in simulations
of correlation structures equals the bias estimated from the genetic
data as in \citet{Yang10}. Since the bias increases with $\epsilon$ (as
the signal-to-noise ratio decreases), finding the appropriate value of
$\epsilon$ should be straightforward using simulations.

The example provided on estimating heritability shows that a
sensitivity analysis can uncover obstacles in applying methods---which
work for simulated data---to actual data collected with error. In the
following section we provide an additional sensitivity analysis applied
to a clustering algorithm.\looseness=1

\subsection{Clustering}

As an additional example, we consider the PAM algorithm [\citet
{Kaufman1990}] to cluster data which has known structure but different
levels of noise. We use the adjusted Rand statistic [\citet{rand},
\citet{Yeung}]
to measure the degree of concordance between the clustering output and
the truth. Using silhouette width, the unsupervised PAM algorithm will
give the optimal number of clusters. The adjusted Rand statistic models
the degree of concordance between the PAM results and the truth. An
adjusted Rand of 1 indicates perfect concordance; an adjusted Rand of
zero indicates a random partition. For each of the models we tested, we
created the correlation matrix (including noise) using an appropriately
customized algorithm.

\subsubsection*{Clustering results from simulation}
\begin{figure}

\includegraphics{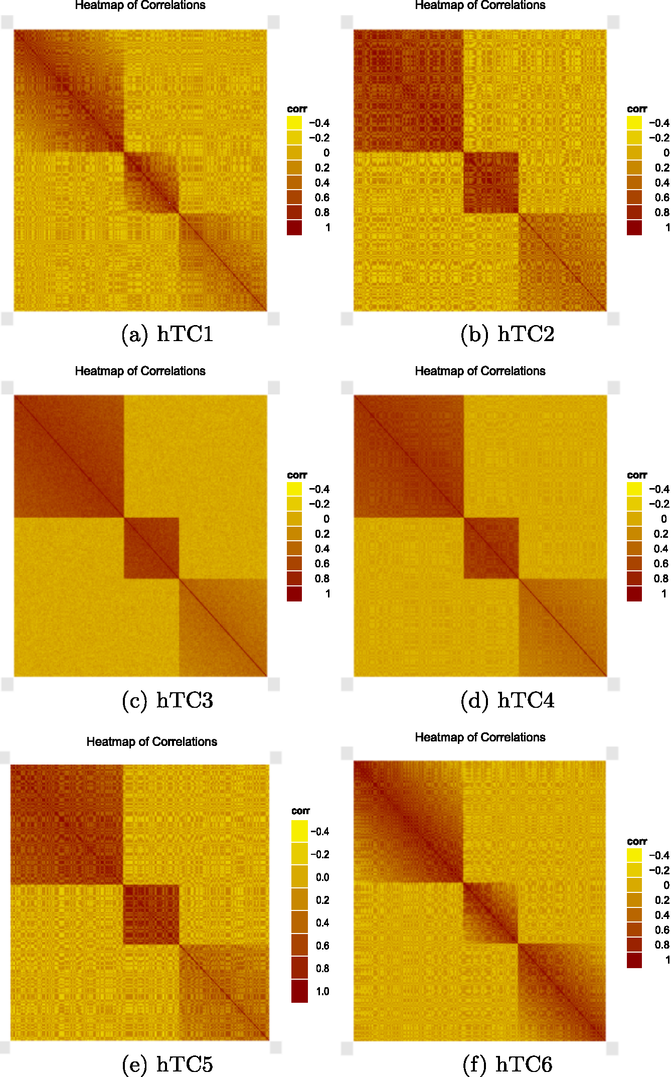}

\caption{Each heatmap represents the
correlation matrix from the scenarios given above.}\label{hTCheatmap}
\end{figure}

Using the hub correlation structure discussed in Section \ref
{subsubhub}, we are required to choose a method to fill out the rest of
the correlation matrix. We use the Toeplitz structure as discussed in
Algorithm \ref{AlgorithmToeplitzHub} with the parameter settings below.
All simulations were done in the three cluster setting with groups of
the following size: $g_1=100, g_2=50, g_3=80$. Recall that with the
hub-Toeplitz correlation, the correlation values descend according to
some power (here linearly) from a specified maximum to a specified
minimum correlation (see Figure \ref{hTCheatmap}):\vspace*{9pt}

\hspace*{-14pt}\noindent\tabcolsep=0pt
{\fontsize{10.5}{12.5}{\selectfont{
\begin{tabular*}{\textwidth}{@{\extracolsep{\fill}}lccccc@{}}
\hline
\mbox{(a) hTC1} &$\rho_1
\in(0.7 \rightarrow0)$ & $\rho_2 \in(0.7 \rightarrow0)$ &
$\rho_3 \in(0.4 \rightarrow0)$ & $\mathbf{u}_i \in
\R^{2}$ & $\epsilon= 0.23$
\\[3pt]
\mbox{(b) hTC2} &$\rho_1 \in(0.7 \rightarrow0.5)$&
$\rho_2 \in(0.7 \rightarrow0.6)$& $\rho_3 \in(0.4
\rightarrow0.2)$& $\mathbf{u}_i \in\R ^{2}$& $\epsilon= 0.29$
\\[3pt]
\mbox{(c) hTC3} &$\rho_1 \in(0.7 \rightarrow0.5)$&
$\rho_2 \in(0.7 \rightarrow0.6)$& $\rho_3 \in(0.4
\rightarrow0.2)$& $\mathbf{u}_i \in\R ^{25}$& $\epsilon= 0.29$
\\[3pt]
\mbox{(d) hTC4} &$\rho_1 \in(0.7 \rightarrow0.5)$&
$\rho_2 \in(0.7 \rightarrow0.6)$& $\rho_3 \in(0.4
\rightarrow0.2)$& $\mathbf{u}_i \in\R ^{2}$& $\epsilon= 0.1$
\\[3pt]
\mbox{(e) hTC5} &$\rho_1 \in(0.7 \rightarrow0.5)$&
$\rho_2 \in(0.7 \rightarrow0.6)$& $\rho_3 \in(0.4
\rightarrow0.2)$& $\mathbf{u}_i \in\R ^{2}$& $\epsilon= 0.25$
\\[3pt]
 \mbox{(f) hTC6} &$\rho_1 \in(0.8 \rightarrow0)$&
$\rho_2 \in(0.75 \rightarrow0)$& $\rho_3 \in(0.7
\rightarrow0)$& $\mathbf{u}_i \in\R^{2}$& $\epsilon= 0.19$
\\[3pt]
\hline
\end{tabular*}}}}\vspace*{9pt}

For each of the scenarios above, we simulated 1000 correlation
matrices. We then clustered the data using PAM; the clustering results
were assessed by determining the number of clusters the algorithm
produced (truth was 3 clusters) as well as the concordance between the
clustering results and the truth (1 gives perfect concordance).

\begin{table}[b]
\caption{Results from
optimal number of clusters as well as the adjusted Rand. The original
correlation structure had 3 clusters. A perfect allocation of points
gives an adjusted Rand of 1.}\label{PAMresultsToeplitzHub}
\begin{tabular*}{\textwidth}{@{\extracolsep{\fill}}ld{2.3}d{2.3}ccd{2.3}c@{}}
\hline
\multicolumn{1}{@{}l}{\textbf{Scenario}} & \multicolumn{1}{c}{\textbf{hTC1}} & \multicolumn{1}{c}{\textbf{hTC2}} &
\multicolumn{1}{c}{\textbf{hTC3}} & \multicolumn{1}{c}{\textbf{hTC4}} &
\multicolumn{1}{c}{\textbf{hTC5}} & \multicolumn{1}{c@{}}{\textbf{hTC6}}\\
\hline
min \# clusters & 3 & 3 & 3 & 3 & 3 & 3\\
median \# clusters & 11 & 8 & 3 & 3 & 3 & 3\\
max \# clusters & 20 & 13 & 3 & 3 & 10 & 3\\
median adj Rand & 0.320 & 0.414 & 1 & 1 & 0.770 & 1\\
\hline
\end{tabular*}
\end{table}

Our results show that adding noise can create scenarios about which the
algorithm is unable to determine the true structure (hTC1 and hTC2) and
scenarios where the noise is not sufficient to decrease the performance
of the algorithm (hTC3, hTC4 and hTC6), as well as situations that work
only sometimes (hTC5) (see Table \ref{PAMresultsToeplitzHub}). For correlation structures that degrade all the
way to zero (hTC1 and hTC6), the algorithm is able to discern the
structure if the original correlations are large (hTC6). For
correlation structures that degrade only a small amount (hTC2, hTC3,
hTC4), the results are based on the amount of error and the dimension
from which the noise vectors from Algorithm \ref{AlgorithmToeplitzHub}
are selected.

\subsubsection*{Clustering results on Fisher's Iris data}
\begin{figure}

\includegraphics{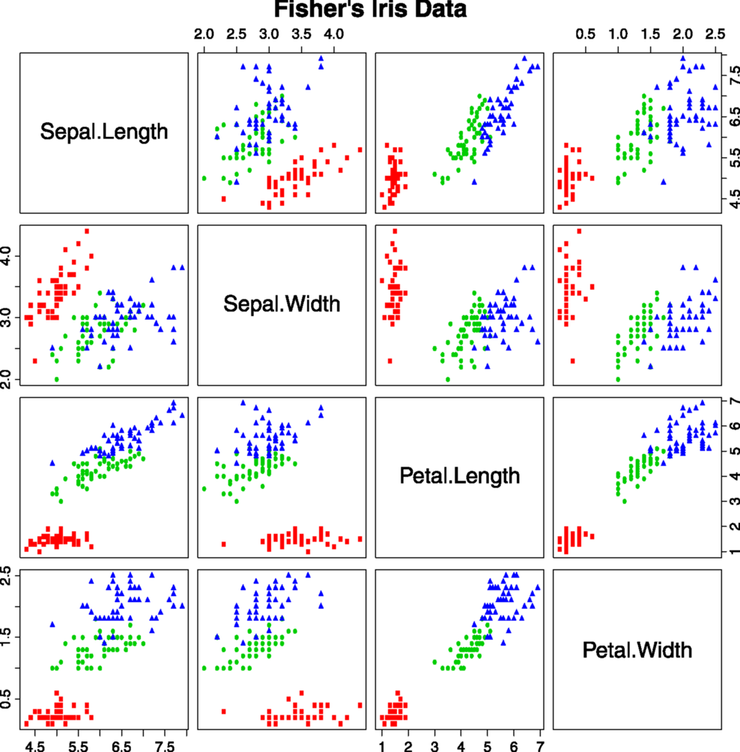}

\caption{The \emph{I. setosa}
flowers are clearly distinct, while \emph{I. versicolor} and \emph{I. virginica} are difficult to differentiate, especially with respect to
sepal measurements.}\label{fisherIrisplot}
\end{figure}

As an application to real data, we consider Fisher's Iris data [\citet
{Fisher36}] which have been used extensively to asses discriminant and
cluster analysis methods. For 50 iris specimens in each of three
species, \emph{Iris setosa}, \emph{I. versicolor} and \emph{I.~virginica},
the sepal length, sepal width, petal length and petal width are
measured in millimeters (see Figure \ref{fisherIrisplot}). Though there are three species measured,
\emph{I.~versicolor} and \emph{I. virginica} are typically quite difficult to
differentiate with unsupervised clustering methods [\citet{Mezzich80},
page 85]. Indeed, when applying the PAM algorithm to the iris data, we
get a perfect separation into two groups (with the three group
silhouette width being slightly smaller).

In order to assess the sensitivity of the output, we can add noise
consistent with the variability of the observations. In this case, we
assume that the correlations within a group are constant, and we
estimate the standard error of the correlations to be approximately
0.01. Such errors might be generated by using vectors from $\R^{25}$,
with $\epsilon= 0.05$. [Refer to the discussion preceding equation
\eqref{SEcor} for choice of $\epsilon$.] Even without noise, the
unsupervised PAM algorithm separates the flowers into two distinct
groups (see Table \ref{PAMresultsIris}). As would be expected with the
iris data, substantially more noise is needed before the two-group
structure is lost.

\begin{table}
\tabcolsep=0pt
\caption{Results for Fisher's Iris
data from optimal number of clusters as well as the adjusted Rand. The
original data have three species, but \emph{I. virginica} and \emph{I.
versicolor} are remarkably similar according to the measurements given
in the data set and are typically grouped together. A perfect
allocation of points gives an adjusted Rand of 1.}\label{PAMresultsIris}
\begin{tabular*}{\textwidth}{@{\extracolsep{\fill}}lccd{2.3}d{2.3}d{2.3}@{}}
\hline
\multicolumn{1}{@{}l}{\textbf{Scenario}} &
\multicolumn{1}{c}{$\bolds{\epsilon=0}$ \textbf{(no error)}} &
\multicolumn{1}{c}{$\bolds{\R^{25}, \epsilon= 0.05}$} &
\multicolumn{1}{c}{$\bolds{ \R^{10}, \epsilon= 0.2}$} &
\multicolumn{1}{c}{$\bolds{\R^{5}, \epsilon= 0.2}$} &
\multicolumn{1}{c@{}}{$\bolds{\R^{2}, \epsilon=0.15}$}\\
\hline
min \# clusters & 2 & 2 & 2 & 2 & 2 \\
median \# clusters & 2 & 2 & 2 & 3 & 4 \\
max \# clusters & 2 & 2 & 3 & 5 & 5\\
median adj Rand&&&&&\\
\quad (2 species) & 1 & 1 & 0.973& 0.570 & 0.469 \\
\hline
\end{tabular*}
\end{table}

\section{Conclusion}\label{SectionConclusion}

We have developed an algorithm for adding noise, in a highly controlled
manner, to a template correlation matrix in order to obtain a more
realistic correlation matrix. Moreover, we have demonstrated how our
general procedure can be tailored to a number of different correlation
models (e.g., constant correlation, Toeplitz structure).

Our method allows for noisy correlation matrices which differ more from
the initial template than the estimated correlation matrix based on
simulated Gaussian data. Using Gaussian data produces a sample
correlation matrix with limited and well-behaved (possibly unrealistic)
differences from the original template correlation if the generated
sample is large. If the generated sample is small, then the sample
correlation matrix is not positive definite (i.e., most of the
eigenvalues will be zero). Using uniform $[-1,1]$ deviates as random
correlation values produces a matrix that is in general not even
positive semidefinite. It can also create relationships between
observations that are meaningless (e.g., a trio of observations where
the first is highly correlated to the other two, but the other two are
negatively correlated).

Using a recent---though already influential---model for heritability
and a standard clustering algorithm, we have shown that simulated
correlation matrices can be used to assess new or existing statistical
methodology. We provide the user with detailed algorithms to use on
several standard clustering structures, as well as a general algorithm
to apply to any correlation matrix for which the smallest eigenvalue
can be reasonably estimated.

\begin{appendix}
\section*{Appendix}\label{Appendix}

\subsection{Preliminaries}\label{SectionPreliminaries}

Recall that if $A$ is a $N \times N$ symmetric matrix, then each of
its eigenvalues is real and, hence,
we may list them in descending order
\[
\lambda_1(A) \geq\lambda_2(A) \geq\cdots\geq
\lambda_N(A)
\]
where each eigenvalue is repeated according to its multiplicity.
According to this convention, $A$ is positive semidefinite if and only
if $\lambda_N(A) \geq0$
and $A$ is positive definite if and only if $\lambda_N(A) > 0$.

The \emph{norm} of a $N \times N$ matrix $A$ is defined to be
%
\begin{equation}
\label{eq-OperatorNorm} \|A\| = \max_{ \|\vec{v} \|= 1} \|A \vec{v} \|
\end{equation}
which equals $\lambda_1(A)$ if $A$ is positive semidefinite. To be more
specific, the expression \eqref{eq-OperatorNorm} is often called the
\emph{operator norm} to distinguish it
from other frequently used matrix norms (e.g., the \emph{Frobenius norm}).
The \emph{condition number} [\citet{HJ}, page 336] of a symmetric
matrix $A$ is defined to be
\[
\kappa(A) = %
\cases{ \bigl\|A^{-1}\bigr\|\|A\|,& \quad $\mbox{if $A$ is
nonsingular},$
\vspace*{2pt}\cr
\infty,& \quad$\mbox{if $A$ is singular}.$} %
\]
In particular, if $A$ is positive semidefinite, then we have
\[
\kappa(A) = %
\cases{ \displaystyle\frac{ \lambda_1(A) }{ \lambda_N(A) }, &\quad $\mbox{if $
\lambda_N(A) > 0$},$
\cr
\infty,&\quad $\mbox{if $\lambda_N(A) =
0$}.$} %
\]
In the following, we let $I_g$ denote the $g \times g$ identity matrix and
$\mathbh{1}_g$ denote the $g \times g$ matrix whose entries are all
equal to $1$.

\subsection{The basic algorithm}\label{SubsectionBasic}
Given an $N \times N$ prototype correlation matrix $\Sigma= (\Sigma
_{ij})_{i,j=1}^N$,
we might wish to add noise to $\Sigma$ in a computationally efficient
way such that the resulting matrix $S$ is also a correlation matrix.
Furthermore, we might also require effective bounds on the condition
number $\kappa(S)$ of $S$ to ensure that $S$ is a suitable candidate
for certain numerical procedures (e.g., matrix inversion). For example,
in the statistical software \texttt{R}, the default tolerance for
detecting linear dependencies in the columns of a matrix is a condition
number $\leq10^{15}$. The following simple procedure accomplishes this task.

\begin{Algorithm}\label{AlgorithmMain}
Let
\begin{longlist}[1.]
\item[1.]$\Sigma$ be a given $N \times N$ correlation matrix,
\item[2.]$0 < \epsilon< \lambda_N(\Sigma)$ ($\epsilon$ is the maximum
noise level),
\item[3.]$M$ be a positive integer (the dimension of the noise space).
\end{longlist}
Select $N$ unit vectors $\vec{u}_1, \vec{u}_2, \ldots, \vec{u}_N$ from
$\R^M$
and form the $M \times N$ matrix $U=(\vec{u}_1 | \vec{u}_2 | \cdots|
\vec{u}_N)$ whose columns
are the $\vec{u}_i$. The $N \times N$ matrix
%
\begin{equation}
\label{eq-SUUI} S = \Sigma+\epsilon\bigl(U^TU - I\bigr)
\end{equation}
is a correlation matrix whose entries satisfy $|S_{ij} - \Sigma_{ij}|
\leq\epsilon$ for $1 \leq i,j\leq N$
and whose condition number $\kappa(S)$ satisfies
%
\begin{equation}
\label{eq-BasicKappa} %
\fbox{$
\displaystyle\kappa(S) \leq\frac{ \lambda_1(\Sigma) + (N-1)\epsilon}{
\lambda_N(\Sigma) - \epsilon}.$ }
\end{equation}
\end{Algorithm}

We might also desire that $\kappa(S) \leq\kappa_{\max}$ for some
fixed $\kappa_{\max}$,
which depends upon the particular requirements of the software being employed.
From \eqref{eq-BasicKappa}, it is easy to see that any $\epsilon> 0$
satisfying the additional constraint
%
\begin{equation}
\label{eq-KappaMax} %
\fbox{$
\displaystyle\epsilon\leq\frac{ \kappa_{\max} \lambda_N(\Sigma) - \lambda
_1(\Sigma
) }{ \kappa_{\max} + (N-1)}$}
\end{equation}
yields an $S$ such that $\kappa(S) \leq\kappa_{\max}$.

\begin{pf*}{Justification of Algorithm \ref{AlgorithmMain}}
Let $E=U^TU$ so that
\[
E = %
\pmatrix{ 1 & \vec{u}_1^T
\vec{u}_2 & \cdots& \vec{u}_1^T
\vec{u}_N\vspace*{2pt}
\cr
\vec{u}_2^T
\vec{u}_1 & 1 & \cdots& \vec{u}_2^T
\vec{u}_N\vspace *{2pt}
\cr
\vdots& \vdots& \ddots& \vdots\vspace*{2pt}
\cr
\vec{u}_N^T \vec{u}_1 &
\vec{u}_N^T \vec{u}_2 & \cdots& 1 }
\]
and note that $E$ is symmetric and positive semidefinite [i.e.,
$\lambda
_N(E) \geq0$].
Moreover, $E$ is positive definite if and only if the $\vec{u}_i$
are linearly independent [\citet{HJ}, Theorem 7.2.10].

Now recall that Ger\v{s}gorin's Disk theorem [\citet{HJ}, Theorem~6.1.1] asserts that
if $A = (A_{ij})_{i,j=1}^N$ is a $N \times N$ matrix, then for each
eigenvalue $\lambda$ of $A$
there exists a corresponding index $i$ such that
\[
| \lambda- A_{ii} | \leq\mathop{\sum_{ j = 1 }}_{ j \neq i}
^N |A_{ij}|.
\]
By Ger\v{s}gorin's theorem and Cauchy--Schwarz, it follows that every
eigenvalue $\lambda$ of $E$ satisfies
\[
|\lambda- 1 | \leq\mathop{\sum_{ j = 1 }}_{ j \neq i}
^N \bigl| \vec {u}_i^T \vec{u}_j\bigr|
\leq(N-1),
\]
whence $0 \leq\lambda_i(E) \leq N$ for $i=1,2,\ldots,N$.\eject

We next define $S$ by \eqref{eq-SUUI}
and observe that $S$ is of the form
%
\begin{equation}
\label{eq-S} S = %
\pmatrix{ 1 & \Sigma_{12} + \epsilon
\vec{u}_1^T \vec{u}_2 & \cdots& \Sigma
_{1N} + \epsilon\vec{u}_1^T
\vec{u}_N\vspace*{2pt}
\cr
\Sigma_{21} + \epsilon
\vec{u}_2^T \vec{u}_1 & 1 & \cdots& \Sigma
_{2N} +\epsilon\vec{u}_2^T
\vec{u}_N\vspace*{2pt}
\cr
\vdots& \vdots& \ddots& \vdots\vspace*{2pt}
\cr
\Sigma_{N1} + \epsilon\vec{u}_N^T
\vec{u}_1 & \Sigma_{N2} + \epsilon\vec{u}_N^T
\vec{u}_2 & \cdots& 1}.
\end{equation}
In particular, $S$ is our original matrix $\Sigma$ with ``noise'' terms
$\epsilon\vec{u}_i^T\vec{u}_j$ of magnitude
at most $\epsilon$ added to the off-diagonal entries. To analyze the
impact of adding this noise,
we require Weyl's Inequalities [\citet{HJ}, Theorem 4.3.1], which
assert that if $A$ and $B$ are
$N \times N$ symmetric matrices, then
%
\begin{equation}
\label{eq-Weyl} \lambda_j(A) + \lambda_N(B) \leq
\lambda_j(A+B) \leq \lambda _j(A) +
\lambda_1(B)
\end{equation}
for $j = 1,2,\ldots,N$. Applying the lower inequality in \eqref
{eq-Weyl} with $j=N$,
$A = \Sigma- \epsilon I_N$ and $B = \epsilon E$, we obtain
\[
0 < \lambda_N(\Sigma) - \epsilon = \lambda_N(\Sigma-
\epsilon I_N) \leq \lambda_N(\Sigma- \epsilon
I_N) + \lambda_N(\epsilon E) \leq \lambda_N(S),
\]
from which we conclude that $S$ is positive definite. Next, we apply
the upper inequality in \eqref{eq-Weyl}
with $j=1$, which yields
\[
\lambda_1(S) \leq \lambda_1(\Sigma- \epsilon I) +
\lambda _1(\epsilon E) \leq \bigl(\lambda_1(\Sigma) -
\epsilon\bigr)+ N \epsilon = \lambda_1(\Sigma) + (N-1)\epsilon.
\]
Putting this all together, we obtain the estimates
\[
0 < \lambda_N(\Sigma) - \epsilon \leq \lambda_N(S) \leq
\lambda_1(S) \leq \lambda_1(\Sigma) + (N-1)\epsilon.
\]
The inequality \eqref{eq-BasicKappa} follows since
$\kappa(S) = \lambda_1(S) / \lambda_N(S)$.
\end{pf*}

There are several arguments which can be made in favor of adding noise
in this manner.
First of all, the procedure described above is easy to implement
numerically, and it can be
rapidly executed. Moreover, it offers a great deal of flexibility
since the dimension $M$
of the ambient space that the vectors $\vec{u}_1, \vec{u}_2, \ldots,
\vec{u}_N$ are drawn from
and the manner in which these vectors are selected is arbitrary and
can be tailored to the particular
application at hand.
Finally, our method is completely general in the sense that any
positive-definite $N \times N$
matrix $E$ having constant diagonal $1$ can be factored as
$E = U^T U$ where $U$ is some matrix whose columns are unit vectors
(e.g., let $U$ be the positive-semidefinite
square root of $E$).
In other words, regardless of the method one employs to produce a
positive-semidefinite matrix $E = U^T U$ for use
in \eqref{eq-SUUI}, the same $E$ can in principle be generated using
our approach.

Let us now say a few words about the manner in which the vectors $\vec
{u}_i$ are selected. If $M$ is very small (e.g., $2 \leq M \leq5$),
then many of the dot products $\vec{u}_i^T \vec{u}_j$ will be large in
magnitude. For many purposes, this yields a very noisy coefficient
matrix $S$ based upon the original template $\Sigma$. Moreover, even if
$M$ is relatively large, then the matrix $E = U^TU$ can be computed
extremely rapidly since generating the unit vectors $\vec{u}_i$ and
computing the dot products $\vec{u}_i^T \vec{u}_j$ involve
straightforward computations (e.g., no eigenvalue calculations).

There are of course many other ways which one could select the $\vec
{u}_i$. If one wishes
the $\vec{u}_i^T \vec{u}_j$ to be consistently large in magnitude
while also ensuring that $E$ has full rank, one lets $M \geq N$
and then selects numbers $\alpha_1,\alpha_2,\ldots,\alpha_N$ at random
from $[-1,1]$
using a continuous probability density function $f(x)$ on $[-1,1]$
which favors extreme values
(e.g., $f(x) = |x|$, $f(x) = \frac{2 - 2 \sqrt{1 - x^2} }{4 - \pi}$ or
a Beta distribution transformed to exist on the range $[-1,1]$). One then
replaces the numbers $\vec{u}_i^T \vec{u}_j$ in \eqref{eq-S} by
%
\begin{equation}
\label{eq-SquareRoot} \alpha_i \alpha_j + \sqrt{
\bigl(1- |\alpha_i|^2\bigr) \bigl(1 - |
\alpha_j|^2\bigr)} \vec {u}_i^T
\vec{u}_j.
\end{equation}
In effect, one is replacing the $\vec{u}_i \in\R^M$ with the unit vectors
$(\alpha_i, \sqrt{1- |\alpha_i|^2} \vec{u}_i) \in\R^{M+1}$. These
vectors tend to have high negative or positive correlations (but they
are linearly independent) since the numbers $\alpha_i$ favor extreme
values in the interval $[-1,1]$.

\subsection{\texorpdfstring{Justification of Algorithm \protect\ref{AlgorithmBlocks}}
{Justification of Algorithm 1}}\label{SubsectionAlgorithmBlocks}

In order to introduce a significant amount of noise to the
off-diagonal blocks, we work instead
with the modified correlation matrix
%
\begin{equation}
\label{eq-SigmaBlocks} \Sigma' = \underbrace{ %
\pmatrix{
\Sigma_1 - \delta\mathbh{1}_{g_1} &&&\vspace*{2pt}
\cr
&
\Sigma_2 - \delta\mathbh{1}_{g_2}&&\vspace*{2pt}
\cr
&&\ddots&
\vspace*{2pt}
\cr
&&&\Sigma_K - \delta\mathbh{1}_{g_K}}
}_{A} + \delta\mathbh{1}_N
\end{equation}
where $\mathbh{1}_g$ denotes the $g \times g$ matrix whose entries are
all $1$. Since
\[
\Sigma_k - \delta\mathbh{1}_{g_k} = ( 1 -
\rho_k) I_{g_k} + (\rho _k - \delta)
\mathbh{1}_{g_k},
\]
it follows that
%
\begin{equation}
\label{eq-BlockMax} \lambda_j( \Sigma_k - \delta
\mathbh{1}_{g_k}) = %
\cases{ g_k(
\rho_k - \delta) + (1 - \rho_k), & \quad$\mbox{if $j =1$},$
\vspace*{2pt}
\cr
1 - \rho_k, & \quad$\mbox{if $j = 2,3,\ldots,
g_k$},$} %
\end{equation}
and that the eigenspace corresponding to the largest eigenvalue of
$\Sigma_k - \delta\mathbh{1}_{g_k}$ is spanned
by the vector $\mathbf{1}_{g_k} = (1,1,\ldots,1) \in\R^{g_k}$. In
particular, the eigenspace corresponding to the
eigenvalue $1 -\rho_k$ is $(g_k - 1)$-dimensional and any eigenvector
$\vec{v} = (v_1,v_2,\ldots,v_{g_k})$
belonging to this eigenspace is orthogonal to $\mathbf{1}_{g_k}$ (i.e.,
satisfies $\sum_{i=1}^{g_k} v_i = 0$).

If we augment $\vec{v}$ by placing $N-g_k$ zeros appropriately, we
obtain a vector
\[
\vec{v}' = (\underbrace{0,0,\ldots,0}_{g_1+\cdots+g_{k-1}},
v_1,v_2,\ldots,v_{g_k}, \underbrace{0,0,
\ldots,0}_{g_{k+1} + \cdots+ g_K}) \in\R^N
\]
which is an eigenvector of $\Sigma'$ corresponding to the eigenvalue $1
- \rho_k$ since
$A \vec{v}' = (1 - \rho_k) \vec{v}'$ and $\mathbh{1}_N \vec{v}' =
\vec{0}$.
It follows that the lowest $N-K$ eigenvalues of $\Sigma$
are the numbers $1 - \rho_k$, each repeated $g_k - 1$ times. In particular,
\[
\lambda_N\bigl(\Sigma'\bigr) = 1 -
\rho_{\max}.
\]
An upper bound on the eigenvalues of $\Sigma$ follows from \eqref{eq-Weyl}
and \eqref{eq-BlockMax}:
\begin{eqnarray*}
\lambda_1\bigl(\Sigma'\bigr) &\leq&
\lambda_1(A) + \lambda_1(\delta\mathbh{1}_N)
\\
&\leq&\max_{1\leq k \leq K}\bigl\{g_k(\rho_k -
\delta) + (1 - \rho_k) \bigr\} + N \delta
\\
&\leq& N(1-\delta) + 1 + N \delta
\\
&=& N+1.
\end{eqnarray*}
Plugging the matrix $\Sigma'$ into Algorithm \ref{AlgorithmMain} and
using the preceding
estimates for $\lambda_1(\Sigma')$ and $\lambda_N(\Sigma')$ into
\eqref{eq-BasicKappa}, we obtain the desired estimate \eqref
{eq-KappaBlocks} for $\kappa(S)$.\qed

\subsection{\texorpdfstring{Justification of Algorithm \protect\ref{AlgorithmToeplitz}}
{Justification of Algorithm 2}}\label{SubsectionAlgorithmToeplitz}

Using the spectral theory of self-adjoint Toeplitz operators,
it is possible to show that $T_g$ is positive definite and that its eigenvalues
satisfy
%
\begin{equation}
\label{eq-ToeplitzSharp} \frac{1-\rho}{1+\rho} \leq \lambda_j(T_g)
\leq \frac{1+
\rho
}{1-\rho}
\end{equation}
for $j=1,2,\ldots,g$. We also remark that the preceding bounds are
quite sharp in the sense that
%
\begin{equation}
\label{eq-ToeplitzLimit} \lim_{g\to\infty} \lambda_1(T_g)
= \frac{1+\rho}{1-\rho}, \qquad \lim_{g\to\infty} \lambda_g(T_g)
= \frac{1-\rho}{1+\rho}
\end{equation}
as the size $g$ of the matrix tends to infinity. In light of the
explicit bounds \eqref{eq-ToeplitzSharp},
a~straightforward application of Algorithm \ref{AlgorithmMain} yields
the following procedure.

\begin{figure}

\includegraphics{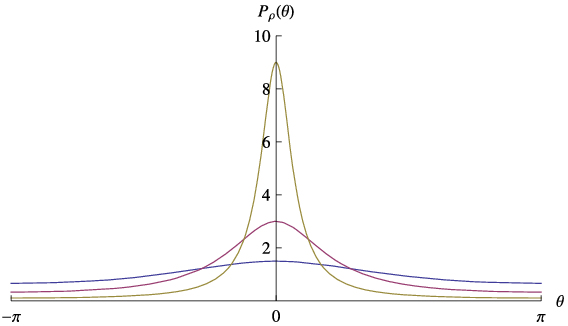}

\caption{The Poisson kernel $P_{\rho}(\theta)$ for $\rho
= 0.2,  0.5,  0.8$. As $\rho\to1^-$,
the graphs spike sharply at $\theta= 0$ while tending rapidly to zero
for $\theta$ away from $0$. Intuitively,
the functions $P_{\rho}(\theta)$ approximate a point mass (i.e., Dirac
$\delta$-function) at $\theta= 0$ as $\rho\to1^-$.}
\label{FigurePoisson}
\end{figure}

To justify the crucial inequalities
\eqref{eq-ToeplitzSharp} and the limits \eqref{eq-ToeplitzLimit},
first observe that the Toeplitz matrix
%
\begin{equation}
\label{eq-Toeplitz} T_g = %
\pmatrix{ 1 & \rho&
\rho^2 & \rho^3 & \cdots& \rho^{g-1}\vspace*{2pt}
\cr
\rho& 1 & \rho& \rho^2 & \cdots& \rho^{g-2}\vspace*{2pt}
\cr
\rho^2 & \rho& 1 & \rho& \cdots& \rho^{g-3}\vspace*{2pt}
\cr
\rho^3 & \rho^2 & \rho& 1 & \cdots&
\rho^{g-4}\vspace*{2pt}
\cr
\vdots& \vdots& \vdots& \vdots& \ddots& \vdots
\vspace*{2pt}
\cr
\rho^{g-1} & \rho^{g-2} & \rho^{g-3} &
\rho^{g-4} & \cdots& 1} %
\end{equation}
is simply the upper-left corner of
the infinite Toeplitz matrix
%
\begin{equation}
\label{eq-BigToeplitz} T = %
\pmatrix{ 1 & \rho& \rho^2 &
\rho^3 & \cdots\vspace*{2pt}
\cr
\rho& 1 & \rho& \rho^2 &
\cdots\vspace*{2pt}
\cr
\rho^2 & \rho& 1 & \rho& \cdots\vspace*{2pt}
\cr
\rho^3 & \rho^2 & \rho& 1 & \cdots\vspace*{2pt}
\cr
\vdots& \vdots& \vdots& \vdots& \ddots} %
\end{equation}
which induces a linear operator $T$ on the Hilbert space $\ell^2$
of all square-summable infinite sequences. Since
the $ij$th entry of $T$ is the $(i-j)$th complex Fourier coefficient
of the function
$P_{\rho}(\theta)\dvtx [-\pi,\pi]\to\R$ defined by
\[
P_{\rho}(\theta) = \sum_{n=-\infty}^{\infty}
\rho^{|n|} e^{in\theta} =\frac{1-\rho^2}{1 - \rho\cos\theta+ \rho^2},
\]
we conclude from \citet{LTTM}, Theorem 1.9, that $T$ is a bounded
self-adjoint operator
whose spectrum equals the range of $P_{\rho}$ [\citet{Halmos}, Problem~250]
[note that $P_{\rho}(\theta)$ is the so-called \emph{Poisson kernel}
from the study of harmonic functions].
A short calculus exercise reveals that $P_{\rho}(\theta)$ achieves
its maximum
value $\frac{1+\rho}{1-\rho}$ at $\theta= 0$ and its minimum value
$\frac{1-\rho}{1+\rho}$
at $\theta= \pm\pi$ (see Figure~\ref{FigurePoisson}),
%
from which we conclude that the spectrum of $T$ is precisely the closed interval
$[ \tfrac{1-\rho}{1+\rho}, \tfrac{1+\rho}{1-\rho} ]$.
By \citet{LTTM}, Proposition 2.17, it follows that the eigenvalues
of $T_n$ are also contained in this interval. This establishes the
inequalities \eqref{eq-ToeplitzSharp}.
The limiting behavior \eqref{eq-ToeplitzLimit} follows immediately
from \citet{LTTM}, Theorem~5.14.\qed

\subsection{\texorpdfstring{Justification of Algorithm \protect\ref{AlgorithmToeplitzHub}}
{Justification of Algorithm 3}}\label{SubsectionAlgorithmToeplitzHub}
By Ger\v{s}gorin's Disk theorem \citet{HJ}, Theorem 6.11, the largest
eigenvalue $\lambda_1(\Sigma_k)$
of $\Sigma_k$ satisfies
\begin{eqnarray*}
\lambda_1(\Sigma_k) &\leq&1 + \rho_k + (
\rho_k - \tau_k) + \cdots+ \bigl(\rho_k -
(g_k-2)\tau _k\bigr)
\\
&=& 1 + (g_k-1)\rho_k - \tau_k
\frac{(g_k-2)(g_k-1)}{2}.
\end{eqnarray*}
This immediately yields \eqref{eq-EVBig}.
On the other hand, it is possible to show that the smallest eigenvalue
of $\Sigma_k$ satisfies
%
\begin{equation}
\label{eq-Bottcher} \lambda_{g_k}(\Sigma_k) \geq1 -
\rho_k - \tfrac{3}{4}\tau_k.
\end{equation}
To be brief, one regards the original $g_k \times g_k$ Toeplitz matrix
$\Sigma_k$ as the
upper-left principal submatrix of a $(2g_k-1) \times(2g_k-1)$
symmetric circulant matrix, the eigenvalues
of which can be exactly computed using well-known techniques [\citet
{MR2179973}, page 32]. A series of
elementary but tedious
algebraic manipulations and a standard eigenvalue interlacing result
[\citet{MR2179973}, Theorem 9.19] yield
the desired inequality \eqref{eq-Bottcher}, from which~\eqref
{eq-EVSmall} follows.
We thank A.~B\"ottcher, the author of \citet{MR2179973},
\citet{LTTM}, for
suggesting this approach to us.
\end{appendix}


\begin{supplement}[id=suppA]
\stitle{R code}
\slink[doi]{10.1214/13-AOAS638SUPP} 
\sdatatype{.r}
\sfilename{aoas638\_supp\_rcode.r}
\sdescription{R code for functions available at
\url{http://pages.pomona.edu/\textasciitilde jsh04747/research/simcor.r}.}
\end{supplement}

%
%

\printaddresses

\end{document}